\documentclass[12pt]{article}

\usepackage{amsmath,amssymb,amsthm,xcolor,hyperref}

\title{CONVOLVED NUMBERS OF K-SECTION OF THE FIBONACCI SEQUENCE: PROPERTIES, CONSEQUENCES}
\author{Khamitov V.\thanks{Odessa Polytechnic National University} 
\and Dmytryshyn D.\thanks{(dmitrishin@opu.edu.ua) Odessa Polytechnic National University} 
\and Gray D.\thanks{(dagray@georgiasouthern.edu) Georgia Southern University}
\and Stokolos A.\thanks{(stokolos@georgiasouthern.edu) Georgia Southern University} 
}

\begin{document}
\date{}
\maketitle

\begin{abstract}
 One possible data encryption scheme is related to stream ciphers, which use a sufficiently long pseudo-random sequence. To increase the cryptographic strength of the cipher, linear shift algorithms are additionally used. Such shifts are generated by linear recurrent sequences. Among them, the most popular are the Fibonacci sequence $\{F_n\}_{n=1}^\infty$  and its generalizations, for example, convolved Fibonacci numbers  $\{F_n^{(s)}\}_{n=1}^\infty$ and k-sections of the  Fibonacci sequence $\{\Phi_{n,k}\}_{n=1}^\infty$ $( \Phi_{n,k}=F_{nk}/F_k).$  This article considers a further generalization of Fibonacci numbers, namely convolutions of k-sections of the Fibonacci sequence $\{\Phi_{n,k}^{(s)}\}_{n=1}^\infty$. 
 These numbers are defined by the relations: 
 $$ 
\Phi_{n,k}^{(1)}=\sum_{j=0}^{n-1}\Phi_{j+1,k}\Phi_{n-j,k\,},\qquad \Phi_{n,k}^{(s)}=\sum_{j=0}^{n-1}\Phi_{j+1,k}\Phi_{n-j,k}^{(s-1)}\,,\quad s=2,3,...
$$
Moreover, $\Phi_{n,1}=F_n, \Phi_{n,1}^{(s)}=F_n^{(s)}$. An explicit formula for the representation of convolvutions of k-sections of the Fibonacci sequence is established:
$$
\Phi_{n,k}^{(s)}=5^{-s}(F_k)^{-2s-1}\sum_{j=0}^{s}(-1)^{(k-1)j}{n+2s\choose j}{n+s-1-j\choose n-1} F_{k(n+2s-2j)},
$$
as well as the Binet type formula: 
\begin{align*}
\Phi_{n,k}^{(s)}=\frac{\varphi^{-k(n+2s)}}{(\varphi^k +\varphi^{-k})^{2s+1}} \sum_{j=0}^{s}&(-1)^{(k-1)j}{n+2s\choose j}{n+s-1-j\choose n-1} \\
& \qquad \cdot (-(-1)^{kn}\varphi^{2kj} + \varphi^{2k(n+2s-j)}),
\end{align*}
$\varphi=\frac12(1+\sqrt 5).$ Several consequences were also obtained for $F_n$  and $F_n^{(s)}$.  These formulas are based on the connection between the derivatives of Chebyshev polynomials of the second kind $U_n^{(s)}(z)$  and the Chebyshev polynomials themselves, as well as the connection for convolutions of k-sections of the  Fibonacci sequence with derivatives of Chebyshev polynomials of the second kind via Lucas numbers $L_k$, namely, 
\[
\Phi_{n,k}^{(s)}=\frac1{2^s s!}  \left\{\begin{array}{@{}l@{}}
   (-i)^{n-1} U^{(s)}_{n+s-1}\left(\frac i2 L_k\right),\qquad k -odd,\\
    U^{(s)}_{n+s-1}\left(\frac 12 L_k\right),\qquad\qquad\quad\;\, k -even
  \end{array}\right.\,.
\]
Note that the sequences $\{\Phi_{n,k}^{(s)}\}_{n=1}^\infty$ for $k=3,4,...$ $s=1,2,..$ are not included in the OEIS encyclopedia.\\


\end{abstract}

\section{Motivation}
One of the main problems of digital data processing is the problem of unauthorized access to information. To protect data in the digital environment cryptography methods are used. Cryptography is an integral part of modern society (electronic mail, online shopping, health information, block-chain technologies  etc.).

There are different data encryption schemes; one of the possible schemes is related to stream ciphers, which use a sufficiently long pseudo-random sequence. A stream scheme converts a stream of text characters into a stream of cipher-text, and the conversion depends on the state of the system. Identical text characters will be encrypted into different ciphertext characters. The simplest stream scheme is the Vernam cipher \cite{ref1}.

To increase the cryptographic strength of the cipher, linear or nonlinear shift algorithms are additionally used. Linear shift gives speed while nonlinear increases cryptographic strength. If a pseudo-random sequence is generated using nonlinear maps, then it is possible to additionally use a linear shift. Such shifts are generated by linear recurrent sequences. Among them, the most popular are the Fibonacci sequence and its generalizations \cite{ref2}. On various aspects of cryptography there exist numerous literatures, for example, \cite{ref3, ref4}. 

	The connection between combinatorial methods inherent in the properties of recurrent sequences and the problem of cryptanalysis once again demonstrates the philosophy of intertwining different branches of mathematics, how from one fascinating theoretical problem solutions to other quite practical problems can be generated. 

\section{Introduction and problem statement}

  	The Fibonacci sequence $\{F_n\}_{n=1}^\infty$ is defined by the recurrence relations 
\begin{equation} \label{eq1} 
x_{n+2}=x_{n+1}+x_{n},     \quad n=1,2...                                %\eqno (1)
\end{equation}
and initial conditions $x_1=1, x_2=1$.

In turn, the Lucas sequence  $\{L_n\}_{n=1}^\infty$ is defined by the same relations but with initial conditions $x_1=1, x_2=3$. On the Online Encylopedia of Integer Sequences (OEIS), the reference numbers for the Fibonacci and Lucas sequences are, respectively, A000045 and A000032.

\begin{table}
\begin{center}
{\renewcommand{\arraystretch}{2}
\begin{tabular}{|l|l|l|l|}
\hline
{\bf Sequence} & {\bf Recursion} & $x_n$, $n = 1,\ldots,10$ & {\bf OEIS \#} \\ \hline
Fibonacci &
\parbox{3cm}{$x_{n + 2} = x_{n + 1} + x_n$ \\
$x_1 = 1$, $x_2 = 1$} &
\parbox{5cm}{
\vspace{0.1cm}

1, 1, 2, 3, 5, 8, 13, 21, 34, 55, \ldots

\vspace{0.1cm}
} &
A000045 \\ \hline
Lucas &
\parbox{3cm}{$x_{n + 2} = x_{n + 1} + x_n$ \\
$x_1 = 1$, $x_2 = 3$} &
\parbox{5cm}{
\vspace{0.1cm}

1, 3, 4, 7, 11, 18, 29, 47, 76, 123, \ldots

\vspace{0.1cm}
} &
A000032 \\ \hline
$F_{2n}$ &
\parbox{4cm}{$x_{n + 2} = 3x_{n + 1} - x_{n}$\\
$x_1 = 1$, $x_2 = 3$} &
\parbox{5cm}{
\vspace{0.1cm}

1, 3, 8, 21, 55, 144, 377, 987, 2584, 6765, \ldots

\vspace{0.1cm}
} &
A001906 \\ \hline
{$F_{2n-1}$} &
\parbox{4cm}{{ $x_{n + 2} = 3x_{n + 1} - x_{n}$\\
$x_1 = 1$, $x_2 = 2$}} &
\parbox{5cm}{
\vspace{0.1cm}

{1, 2, 5, 13, 34, 89, 233, 610, 1597, 4181, \ldots}

\vspace{0.1cm}
} &
{A001519} \\ \hline
3-section &
\parbox{3.5cm}{$x_{n + 2} = 4x_{n + 1} + x_n$\\
$x_1=1$, $x_2=4$} &
\parbox{5cm}{
\vspace{0.1cm}

1, 4, 17, 72, 305, 1292, 5473, 23184, 98209, 416020, \ldots

\vspace{0.1cm}
} &
A001076 \\ \hline
\end{tabular}}
\end{center}
\caption{Various sequences relating to the Fibonacci numbers, from the Online Encyclopedia of Integer Sequences (OEIS) \cite{ref5}.}\label{tab1}
\end{table}

% it might be better to include the various sequences in and their OEIS reference numbers in a common table that we can reference in the text. I can set that up when I get to a computer this afternoon. I am editing from my phone now, so I can only check the English and make minor changes at the moment.

	%According to the OEIS electronic encyclopedia [5] these sequences have their own numbers: for the Fibonacci sequence $\{F_n\}_{n=1}^\infty=\{1,1,2,3,5,8,13,21,34,55,..\},$  the number is A000045;  for the Lucas sequence $\{L_n\}_{n=1}^\infty=\{1,3,4,7,11,18,29,47,,76,,123,..\}$    it is A000032  respectively.

Other examples of linear recurrence sequences of the second order with different initial conditions can be found in \cite{ref6}. For example, the recurrence sequence 
$$ x_{n+2}=3x_{n+1}-x_{n},     \quad n=1,2...                     
$$
with the initial conditions $x_1=1$ and $x_2=3$, defines a section of the  Fibonacci sequence 
 $\{F_{2n}\}_{n=1}^\infty=\{1,3,8,21,55,..\},$   with OEIS  number A001906, while with the initial conditions $x_1=1$ and $x_2=2,$ defines the  section   with OEIS  number  A001519. 
 
Let us turn our attention to another example which will play a critical role in this paper. Denote $ \Phi_{n,k}=F_{nk}/F_k,\; n=1,2,...,\; k=1,2,..$. Since the number $nk$  is a multiple of $k$, then all numbers  $ \Phi_{n,k}$ are integers \cite{ref7}. In addition, for each $k=1,2,...$ the sequence $ \{\Phi_{n,k}\}_{n=1}^\infty$ is defined by the recurrence relation  
\begin{equation}\label{eq2}
x_{n+2}=L_kx_{n+1}-(-1)^kx_{n},     \quad n=1,2...                
\end{equation}                               
and initial conditions $x_1=1, x_2=L_k$,  \cite{ref7,ref8}.

The numbers $\Phi_{n,k}$  can also be expressed directly via the Lucas  numbers  \cite{ref9}
\begin{equation} \label{eq3}
\Phi_{n,k}=\sum_{j=0}^{[\frac{n-1}2]} (-1)^{(k-1)j}{n-1-j\choose j}(L_k)^{n-1-2j}.                    
\end{equation}                                 
For each  $k$, the sequence $\{\Phi_{n,k}\}_{n=1}^\infty$  has its own number in \cite{ref5}, however, it does not have a special name (at least the authors did not find this name in the literature). In the present article, we will call this sequence the $k$-section of the Fibonacci sequence. Corresponding numbers in   OEIS are: $k=1$ , A000045 (Fibonacci sequence), $k=2$ , A001906 (section of the Fibonacci sequence); $k=3$ , A001076; $k=4$ , A004187; $k=5$ , A049666; $k=6$ , A049660; $k=7$  , A049667; $k=8$ , A049668;  $k=9$, A049669, etc. 
%{\color{red} The sequence $A001519$ can also be thought of as the $2$-section of the Fibonacci sequence provided we start at $F_0 = 0$, $F_1 = 1$, and so we include this with other select $k$-sections in Table~\ref{tab1}.}

Note that numbers generated by general second-order recurrence relations
$x_{n+2}=\alpha x_{n+1}\pm x_n$, $n=1,2,...$ are studied sufficiently enough, but the sequences $\{\Phi_{n,k}\}_{n=1}^\infty$ are special in that they are generated by ``the golden 
ratio" (i.e. expressed via numbers  $\varphi=\frac12(1+\sqrt{5})$,  $\varphi^{-1}=\frac12(-1+\sqrt{5})$).

The next generalization of the Fibonacci sequence that is closely linked to our paper is the convolved Fibonacci numbers, which are defined as follows:
\begin{equation}\label{eq4}
F_n^{(0)}=F_n,\qquad F_n^{(s)}=\sum_{j=0}^{n-1}F_{j+1} F_{n-j}^{(s-1)},\qquad s=1,2,3,...
\end{equation}
From formulas \eqref{eq4} it follows that  $F_1^{(s)}=1, F_2^{(s)}=s+1$ . 
Note that in \cite{ref10} the authors use a different indexation of the sequence
$\{F_n^{(s)}\}_{n=1}^\infty.$ Namely, \cite[(1.1)]{ref10} and \cite[(1.2)]{ref10}
imply that $F_n^{(s)}=0$ for $n\le s,$ $F_{s+1}^{(s)}=1,$ and $F_{s+2}^{(s)}=s+1.$ Thus, the sequence from \cite{ref10} coinsides with the sequence defined by \eqref{eq4} with index shifted by $s.$ It seems to us that the form of defining convolved numbers via \eqref{eq4} is more convenient, and such a form will be used in the article below. 

%{\color{blue}  Note that in \cite{ref10} the author uses indexation starting from zero and ending at $n$. {\color{red} Accounting for a slight shift in indices (\cite{ref10} has $F_0 = 1$ and $F_1 = 1$ while A000045 has $F_0 = 0$ and $F_1 = 1$)} it seems to us that the form of defining convolved numbers via \eqref{eq4} is more convenient, and such a form will be used in the article below. }

\begin{table}
\begin{center}
{\renewcommand{\arraystretch}{1.5}
\begin{tabular}{|l|l|l|}
\hline
{\bf s} & {$F_n^{(s)}$, $n = 1,\ldots, 10$} & {\bf OEIS \#}\\ \hline
1 & 1, 2, 5, 10, 20, 38, 71, 130, 235, 420, \ldots & A001629 \\ \hline
2 & 1, 3, 9, 22, 51, 111, 233, 474, 942, 1836, \ldots & A001628 \\ \hline
3 & 1, 4, 14, 40, 105, 256, 594, 1324, 2860, 6020, \ldots & A001872 \\ \hline
4 & 1, 5, 20, 65, 190, 511, 1295, 3130, 7285, 16435, \ldots & A001873 \\ \hline
\end{tabular}}
\end{center}
\caption{The convolved Fibonacci numbers for $s = 1$, $s = 2$, $s = 3$, and $s = 4$.}\label{tab2}
\end{table}

 For each $s$, the sequence $\{F_{n}^{(s)}\}_{n=1}^\infty$  has its own OEIS number (see Table~\ref{tab2}).

 % $\{F_{n}^{(1)}\}_{N=1}^\infty=\{1,2,5,10,20,38,...\}$   - A001629;   $\{F_{n}^{(2)}\}_{N=1}^\infty=\{1,3,9,22,51,111,...\}$ - A001628;  
 
 % $\{F_{n}^{(3)}\}_{N=1}^\infty=\{1,4,14,40,105,256,...\}$  - A001872;   $\{F_{n}^{(4)}\}_{N=1}^\infty=\{1,5,20,65,190,511,...\}$ - A001873 etc. 

Convolved Fibonacci numbers were studied in several papers - see some references, for example, in \cite{ref5}. The sequence  $\{F_{n}^{(s)}\}_{n=1}^\infty$  is defined by a linear recurrence relation of order higher than the second (for most $s$  the order is equal to $2s$). For example, for $s=1$ we see the convolved numbers satisfy the relation \cite{ref10}:
$x_{n+4}=2x_{n+3}+x_{n+2}-2x_{n+1}-x_n,\; n=1,2,...$
These recurrence relations are easy to obtain by the method of generating functions (this method will be considered further). The connection between generating functions and recurrence relations is described in detail, for example, in \cite{ref7,ref11}.

For $k$-sections of the Fibonacci sequence it is also possible to define convolved numbers by a formula similar to (4), namely
\begin{equation}\label{eq5}
\Phi_{n,k}^{(1)}=\sum_{j=0}^{n-1}\Phi_{j+1,k} \Phi_{n-j,k}\; ,\quad
\Phi_{n,k}^{(s)}=\sum_{j=0}^{n-1}\Phi_{j+1,k} \Phi_{n-j,k}^{(s-1)}\; ,\quad s=2,3,...
\end{equation}
Since $\Phi_{n,1}=F_n$, then $\Phi_{n,1}^{(s)}=F_n^{(s)}$. We will explore this sequence in this article.

The main goals of the presented work are:
\begin{itemize}
\item finding representation formulas of convolved $k$-sections of the Fibonacci sequence; 
\item determining properties of the sequence $\{\Phi_{n,k}^{(s)}\}_{n=1}^\infty$ by finding new connections between elements of the sequence;
\item finding of Binet type formulas for the convolved $k$-sections,   $\Phi_{n,k}^{(s)}$.
\end{itemize}

To achieve these goals, the formulas connecting derivatives of the Chebyshev polynomials of the second kind and the Chebyshev polynomials of the second kind themselves were used \cite{ref12}. As side results, various identities related to the Fibonacci and Lucas numbers, binomial coefficients, convolved Fibonacci numbers and k-sections of the Fibonacci sequence were obtained. 

	Additionally, we note that the sequences  $\{\Phi_{n,k}^{(s)}\}_{n=1}^\infty$  for $k=3,4,...$ , $s=1,2,...$   are not represented at the OEIS encyclopedia \cite{ref5}.

\section{Preliminary results}

Relation (1) can be considered as a linear difference equation with the given initial conditions. The roots of the corresponding characteristic equation are $\varphi$  and  $-\varphi^{-1}$, and the particular solution is represented in the form 
\begin{equation}\label{eq6}
F_n=\frac{\varphi^n-(-\varphi)^{-n}}{\varphi+\varphi^{-1}}=\frac1{\sqrt 5}(\varphi^n-(-\varphi)^{-n}).
\end{equation}
Similarly,
\begin{equation}\label{eq7}
L_n=\varphi^n-(-\varphi)^{-n}.
\end{equation}
Formulas \eqref{eq6} and \eqref{eq7} are called Binet formulas for the Fibonacci and Lucas numbers \cite{ref7}.

The function $f(z)=\frac z{1-z-z^2}$  is the generating function for the Fibonacci sequence \cite{ref7}; this means that
 $$f(z)=\frac z{1-z-z^2}=\sum_{j=1}^\infty F_jz^j.
 $$ 
For the section (or $2$-section, in our nomenclature) of the Fibonacci sequence   $\{F_{2n}\}_{n=1}^\infty$  the generating function will be $\displaystyle\hat f(z)=\frac z{1-3z+z^2}$, where  $\displaystyle \frac z{1-3z+z^2}=\sum_{j=1}^\infty F_{2j}z^j.$

	Now consider the $k$-section of the Fibonacci sequence. Using \eqref{eq2} and the connection between the generating function and recurrence relations, one can determine the generating function for $\Phi_{j,k}$  \cite[p.~230]{ref7}:  $\displaystyle\hat f_k(z)=\frac z{1-L_kz+(-1)^kz^2}$, where   
    $$
    \frac z{1-L_kz+(-1)^kz^2}=\sum_{j=1}^\infty \Phi_{j,k}z^j
    $$
( $\hat f_0(z)=f(z),\; \hat f_1(z)=\hat f(z)$).

For the convolved Fibonacci numbers, $F_n^{(s)}$,  the generating function is defined as \cite{ref10} $\displaystyle \tilde f(z)=\frac z{(1-z-z^2)^{s+1}}$  ; $\displaystyle \frac z{(1-z-z^2)^{s+1}}= \sum_{j=1}^\infty F_{j}^{(s)}z^j.$  

The function $g(z,t)=\frac1{1-2tz+z^2}$  is called the generating function for Chebyshev polynomials of the second kind $U_n(t)$  \cite{ref13}, so that
$$
\frac1{1-2tz+z^2}=\sum_{j=0}^\infty U_{j}(t)z^j,
$$ 
where \cite{ref14} 
$$
U_n(t)=\sum_{j=1}^{[\frac n2]}(-1)^j{n-j\choose j}2^{n-2j}t^{n-2j}.
$$
	Since  $zg(z/i,i/2)=f(z),\; zg(z,3/2)=\hat f(z),\; (i^2=-1)$, then \cite{ref7}
\begin{equation}\label{eq8}
(-i)^{n-1} U_{n-1}\left(\frac i2\right)=F_n,
\end{equation}
\begin{equation}\label{eq9}
U_{n-1}\left(\frac 32\right)=F_{2n}. 
\end{equation}                                               
Further, $zg(z/i,(i/2)L_k)=\hat f_k(z),$   if  $k$ is odd, and $zg(z,(1/2)L_k)=\hat f_k(z),$ if $k$ is even. Then, 
\begin{equation}\label{eq10}
\Phi_{n,k}=\frac{F_{nk}}{F_k}=\left\{\begin{array}{@{}l@{}}
   (-i)^{n-1} U_{n-1}\left(\frac i2 L_k\right),\qquad k -odd,\\ \\
    U_{n-1}\left(\frac 12 L_k\right),\qquad\qquad\quad\;\, k -even
  \end{array}\right.\,.
\end{equation}

Let us return to the convolved Fibonacci numbers, $F_{n}^{(s)}$, and let us find 
\begin{align*}
\frac{\partial ^s}{\partial t^s}g(z,t)&=2^ss!\frac{z^s}{(1-2tz+z^2)^{s+1}}=
\sum_{j=0}^\infty U_{j}^{(s)}(t)z^j=\sum_{j=s}^\infty U_{j}^{(s)}(t)z^j\\
&= \sum_{j=1}^\infty U_{j+s-1}^{(s)}(t)z^{j+z-1} \qquad (\mbox{where } U_{j}^{(s)}(t)=0, j<s).
\end{align*}
Then, 
$$
\frac{z}{(1-2tz+z^2)^{s+1}}=\frac1{2^ss!}\sum_{j=1}^\infty U_{j+s-1}^{(s)}(t)z^j,
$$
and consequently, 
\begin{equation}\label{eq11}
F_n^{(s)}=\frac{(-i)^{n-1}}{2^s s!}  U_{n+s-1}^{(s)}\left(\frac i2\right).
\end{equation}
Similarly,
\begin{equation}\label{eq12}
F_{2n}^{(s)}=\frac{(-i)^{n-1}}{2^s s!}  U_{n+s-1}^{(s)}\left(\frac 32\right).
\end{equation}
Formulas \eqref{eq10}, \eqref{eq11}, and \eqref{eq12} are a generalization of the well-known formulas \eqref{eq8} and \eqref{eq9}. 
	
\section{Main results}

Further results  are essentially obtained by using the formulas connected to derivatives of Chebyshev polynomials of the second kind, $U_n^{(s)}(z)$  ($s$ is the order of the derivative), and the Chebyshev polynomials of the second kind themselves, $U_n(z)$. Let us present these formulas \cite{ref12}.
\begin{equation}\label{eq13}
U_{n}^{(s)}(z) =s!\sum_{j=0}^{[\frac{n-s}2]}(-1)^{j}{n-j\choose j}{n-2j\choose s} 2^{n-2j}z^{n-2j-s},
\end{equation}
\begin{equation}\label{eq14}
\frac1{2^ss!}U_{n+s-1}^{(s)}\left(\frac 12(z^{\frac12}+z^{-\frac12})\right)  =
z^{-\frac{n-1}2} \sum_{j=0}^{n-1}{n+s-1-j\choose s}{s+j\choose s} z^j,
\end{equation}
\begin{eqnarray}\label{eq15}
\lefteqn{\frac1{2^ss!}U_{n+s-1}^{(s)}\left(\frac 12(z^{\frac12}+z^{-\frac12})\right)} \nonumber \\ &=&
z^{-\frac{n-1}2}(1-z)^{-(2s+1)} \sum_{j=0}^{s}(-1)^j {n+2s \choose j}{n+s-1-j\choose n-1} \nonumber \\
& & \qquad \qquad \cdot (z^j-z^{n+2s-j}),
\end{eqnarray}
\begin{eqnarray}\label{eq16}
\lefteqn{\frac{2^s}{s!}(1-z^2)^sU_{n + s - 1}^{(s)}(z)} \nonumber \\ 
&=& (-1)^s\sum_{j = 0}^s(-1)^j{n + s - 1 - j \choose n - 1}{n + 2s \choose j}U_{n + 2s - 1 - 2j}(z). 
\end{eqnarray}

Binet type formulas for the convolved Fibonacci numbers are obtained from \eqref{eq11}, \eqref{eq14}, and \eqref{eq15} by substituting $-\varphi^2$ for $z$:
$$
F_{n}^{(s)} = (-1)^{n - 1}\varphi^{-n + 1}\sum_{j = 0}^{n - 1}(-1)^j{n + s - 1 - j \choose s}{s + j \choose s}\varphi^{2j},
$$

\begin{equation}\label{eq17}
F_{n}^{(s)} = \frac{\varphi^{-n + 1}}{(1 + \varphi^2)^{2s + 1}}\sum_{j = 0}^{s}{n + 2s \choose j}{n + s - 1 - j \choose n - 1}\big(-(-1)^n\varphi^{2j} + \varphi^{2(n + 2s - j)}\big). 
\end{equation}
The connection between convolved Fibonacci numbers and the normal Fibonacci number follows from formulas \eqref{eq17} and \eqref{eq6}:

\begin{equation}\label{eq18}
F_n^{(s)} = 5^{-s}\sum_{j = 0}^{s}{n + 2s \choose j}{n + s - 1 - j \choose n - 1}F_{n + 2s - 2j}. 
\end{equation}

In the derivation of \eqref{eq18}, the formulas $\displaystyle \frac{\varphi^{-n + 1}}{(1 + \varphi^2)^{2s + 1}} = \frac{\varphi^{-n}\varphi^{-s}}{(\varphi + \varphi^{-1})^{2s}(\varphi + \varphi^{-1})}$, and $\displaystyle \frac{-(-1)^n\varphi^{-(n + 2s - 2j)} + \varphi^{n + 2s - 2j}}{\varphi + \varphi^{-1}} = F_{n + 2s - 2j}$ were used.

Formula \eqref{eq18} can also be obtained from 
\eqref{eq16} by substituting $z = i/2$. This formula substantially refines formula (1.12) from \cite{ref10}.

For $\displaystyle s = 1$: $F_{n}^{(1)} = \frac{1}{5}(nF_{n + 2} + (n + 2)F_n).$

For $\displaystyle s = 2$: $F_{n}^{(2)} = \frac{1}{25}\left(\frac{n(n + 1)}{2}F_{n + 4} + n(n + 4)F_{n + 2} + \frac{(n + 3)(n + 4)}{2} F_n\right).$

Note that the sequences $\{F_n^{(s)}\}_{s = 1}^{\infty}$ are represented in OEIS by the following numbers for $n = 2$ (A000027), $n = 3$ (A000096) , $n = 4$ (A006503), $n = 5$ (A006504), etc.

Now, let us differentiate identity \eqref{eq16}
$$
\frac{2^s}{s!}s(1 - z^2)^{s - 1}(-2z)U_{n + s - 1}^{(s)}(z) + \frac{2^s}{s!}(1 - z^2)^{s}U_{n + s - 1}^{(s + 1)}(z) = 
$$
$$
(-1)^s\sum_{j = 0}^{s}(-1)^j{n + s - 1 - j \choose n - 1}{n + 2s \choose j}U_{n + 2s - 1 - 2j}^{(1)}(z),
$$
and substitute $z = i/2$. The applying formula \eqref{eq11}, we obtain
$$
4s5^{s - 1}F_n^{(s)} + 2(s + 1)5^sF_{n - 1}^{(s + 1)} = 2\sum_{j = 0}^{s}{n + 2s \choose j}{n + s - 1 - j \choose n -1}F_{n + 2s - 2j - 1}^{(1)},
$$
whence
\begin{eqnarray}\label{eq19}
F_n^{(s)} &=& -\frac{2(s - 1)}{5s}F_{n + 1}^{(s - 1)} \nonumber \\
& & \quad + \frac{1}{5^{s - 1}s}\sum_{j = 0}^{s - 1}{n + 2s - 1 \choose j}{n + s - 1 - j \choose n}F_{n + 2s - 2j - 2}^{(1)}.
\end{eqnarray}
Formula \eqref{eq19} allows us to obtain various corollaries. For example, letting $s = 2$:
$$
F_n^{(2)} = \frac{1}{10}\left((n + 1)F_{n + 2}^{(1)} - 2F_{n + 1}^{(1)} + (n + 3)F_{n}^{(1)}\right).
$$
By taking into account \eqref{eq18}, formula \eqref{eq19} takes the form 
\begin{align*}
F_n^{(s)} = \frac{1}{5^{s - 1}s}\sum_{j = 0}^{s - 1}&{n + 2s - 1 \choose j}{n + s - 1 - j \choose n} \\ & \cdot \left(F_{n + 2s - 2j - 2}^{(1)} - \frac{2(s - 1)}{5}F_{n + 2s -2j - 1}\right).
\end{align*}
Using \eqref{eq19}, it is not difficult to obtain other identities that connect to the convolved Fibonacci numbers.

For $k$-sections of the Fibonacci sequence, the  Binet formuas are obvious. From formula (11) for $s = 0$ and $\displaystyle z = \begin{cases}
(i/2)L_k, & k\mbox{ is odd}\\
(1/2)L_k, & k\mbox{ is even}
\end{cases},$ formula \eqref{eq3} follows immediately. Note that formula \eqref{eq10} is a consequence of formulas \eqref{eq7} and \eqref{eq13} for $s = 0$ and $z = -\varphi^{2k}$.

For convolved $k$-sections of the Fibonacci sequence, $\Phi_{n,k}^{(s)}$, the generating function is defined similarly to the case $F_n^{(s)}$:
$$
f_k^{(s)} = \frac{z}{\left(1 - L_kz + (-1)^kz^2\right)^{s + 1}} = \sum_{j = 1}^{\infty}\Phi_{j,k}^{(s)}z^j.
$$

Indeed, using \cite[p.~216]{ref7} and \cite{ref11}, we obtain
\begin{align*}
\frac{z}{1 - L_kz + (-1)^kz^2} \cdot \frac{1}{1 - L_kz + (-1)^kz
^2} \\ = \sum_{j = 1}^{\infty}\left(\sum_{\ell = 0}^{n - 1} \Phi_{\ell + 1,k}\Phi_{n - \ell,k}\right)z^j = \sum_{j = 1}^{\infty}\Phi_{j,k}^{(1)}z^j,
\end{align*}
\begin{align*}
\frac{z}{(1 - L_kz + (-1)^kz^2)^s} \cdot \frac{1}{1 - L_kz + (-1)^kz
^2} \\ = \sum_{j = 1}^{\infty}\left(\sum_{\ell = 0}^{n - 1} \Phi_{\ell + 1,k}\Phi_{n - \ell,k}^{(s - 1)}\right)z^j = \sum_{j = 1}^{\infty}\Phi_{j,k}^{(s)}z^j.
\end{align*}
Whence,
\begin{equation}\label{eq20}
\Phi_{n,k}^{(s)} = \frac{1}{2^ss!}\begin{cases}
\displaystyle (-i)^{n - 1}U_{n + s - 1}^{(s)}\left(\frac{i}{2}L_k\right), k \mbox{ is odd} \\
\displaystyle U_{n + s - 1}^{(s)}\left(\frac{1}{2}L_k\right), k \mbox{ is even}
\end{cases}. 
\end{equation}

Binet type formulas for the $\Phi_{n,k}^{(s)}$ are obtained from \eqref{eq14}, \eqref{eq15}, and \eqref{eq20} by substituting $z = -\varphi^{2k}$ if $k$ is odd and $z = \varphi^{2k}$ if $k$ is even:
$$
\Phi_{n,k}^{(s)} = (-1)^{k(n - 1)}\varphi^{k(-n + 1)}\sum_{j = 0}^{n - 1}(-1)^{jk}{n + s - 1 - j \choose s}{s + j \choose s}\varphi^{2jk},
$$
\begin{align*}
\Phi_{n,k}^{(s)} = \frac{\varphi^{-k(n + 2s)}}{\left(\varphi^k + \varphi^{-k}\right)^{2s + 1}}\sum_{j = 0}^{s}(-1)^{(k - 1)j}{n + 2s \choose j}{n + s - 1 - j \choose n - 1} \\ \cdot \left(-(-1)^{kn}\varphi^{2kj} + \varphi^{2k(n + 2s - j)}\right).
\end{align*}

From \eqref{eq16} and \eqref{eq6} the following formulas is obtained, which is a complete analogy of formula \eqref{eq18}:
\begin{equation}\label{eq21}
\Phi_{n,k}^{(s)} = 5^{-s}(F_k)^{-2s - 1}\sum_{j = 0}^{s}(-1)^{(k - 1)j}{n + 2s \choose j}{n + s - 1 - j \choose n - 1}F_{k(n + 2s - 2j)}.
\end{equation}

An equivalent formula is:
$$
\Phi_{n,k}^{(s)} = 5^{-s}(F_k)^{-2s}\sum_{j = 0}^{s}(-1)^{(k - 1)j}{n + 2s \choose j}{n + s - 1 - j \choose n - 1}\Phi_{(n + 2s - 2j), k}.
$$

For the particular cases $s = 1$ and $s = 2$, we have:
$$
\Phi_{n,k}^{(1)} = \frac{1}{5(F_k)^2}\left(n\Phi_{n+2,k} + (-1)^{k - 1}(n + 2)\Phi_{n,k}\right);
$$
\begin{align*}
\Phi_{n,k}^{(2)} = \frac{1}{25(F_k)^4}\left(\frac{n(n+1)}{2}\Phi_{n+4,k} + (-1)^{k-1}n(n + 4)\Phi_{n+2,k} \right. \\ \left. + \frac{(n + 3)(n + 4)}{2}\Phi_{n,k}\right).
\end{align*}

Note that $\displaystyle \Phi_{3,1}^{(s)} = \frac{1}{2}(s + 1)(s + 4)$, $\Phi_{3,2}^{(s)} = \frac{1}{2}(s + 1)(9s + 16)$; the sequence $\{\Phi_{3,3}^{(s)}\} = \{50,99,164,245,342,\ldots\}$ counts the "number of walks on a cubic lattice" \cite{ref5}, its OEIS number is A005570.

The sequences $\{\Phi_{n,k}^{(s)}\}_{n = 1}^{\infty}$ for $k = 3, 4, \ldots$, $s = 1, 2, \ldots$ are not represented in the OEIS encyclopedia \cite{ref5}. The sequence
$$
\{\Phi_{n,3}^{(1)}\}_{n = 1}^{\infty} = \{1, 8, 50, 280, 1476, 7472, 36836, \ldots\},
$$
was checked in OEIS, but did not match any known sequences. Similarly, no matches were found for any other values of $s$. Interestingly, not even $\{\Phi_{n,k}\}_{n = 1}^{\infty}$, $k = 13, 14, 15, \ldots$, are represented in OEIS.

It follows from \eqref{eq10}, \eqref{eq13}, and \eqref{eq21} that
$$
\Phi_{n,k}^{(s)} = \sum_{j = 0}^{\left[\frac{n-1}{2}\right]}(-1)^{(k-1)j}{n + s - 1 - j \choose j}{n + s - 1 - 2j \choose s}(L_k)^{n - 1 - 2j}.
$$
Then,

\begin{align}\label{eq22}
\sum_{j = 0}^{\left[\frac{n-1}{2}\right]}(-1)^{(k-1)j}{n + s - 1 - j \choose j}{n + s - 1 - 2j \choose s}(L_k)^{n - 1 - 2j} = 
\nonumber \\
5^{-s}(F_k)^{-2s - 1}\sum_{j = 0}^{s}(-1)^{(k-1)j}{n + 2s \choose j}{n + s - 1 - j \choose n - 1}F_{k(n + 2s - 2j)}.
\end{align}

\begin{table}
\begin{center}
{\renewcommand{\arraystretch}{1.5}
\begin{tabular}{|l|l|}
\hline
{\bf s} & {$\Phi_{n,3}^{(s)}$, $n =1,\ldots, 10$}  \\ \hline
1 & 1, 8, 50, 280, 1475, 7472, 36836, 178000, 847045, 3982200\ldots \\ \hline
2 & 1, 12, 99, 688, 4326, 25464, 143018, 775536, 4092207, 21126564, \ldots\\ \hline
3 & 1, 16, 164, 1360, 9930, 66544, 419124, 2518576, 14585635, 81987680,\ldots \\ \hline
\end{tabular}}
\end{center}
\caption{The convolved $3$-sections of Fibonacci numbers for $s = 1$, $s = 2$, and $s = 3$.}\label{tab3}
\end{table}

Some numbers from $\Phi_{n,3}^{(s)}$ sequences are provided in Table~\ref{tab3}.

\section{Description of the convolution of $k$-sections for OEIS}

Let us summarize our findings for $\{\Phi_{n,k}^{(s)}\}_{n = 1}^{\infty}$. We have the following:
\begin{itemize}
\item by definition
\begin{align*}
\Phi_{n,k} = F_{nk}/F_{k}, \ \Phi_{n,k}^{(1)} = \sum_{j = 0}^{n - 1}\Phi_{j+1,k}\Phi_{n-j,k}, \\ \Phi_{n,k}^{(s)} = \sum_{j = 0}^{n - 1}\Phi_{j+1,k}\Phi_{n-j,k}^{(s-1)}, \ s = 2,3,\ldots;
\end{align*}
\item via generating function
$$
f(z) = \frac{z}{(1 - L_kz + (-1)^kz^2)^{s + 1}} = \sum_{n = 0}^{\infty}\Phi_{n,k}^{(s)}z^n;
$$
\item by Binet formula
\begin{align*}
\Phi_{n,k}^{(s)} &= \frac{\varphi^{-k(n + 2s)}}{(\varphi^k + \varphi^{-k})^{2s + 1}}\sum_{j = 0}^{s}(-1)^{(k-1)j}{n + 2s \choose j}{n + s - 1 - j \choose n - 1} \\
& \qquad \qquad \cdot \left(-(-1)^{kn}\varphi^{2kj} + \varphi^{2k(n + 2s - j)}\right);
\end{align*}
\item the explicit formula via the Fibonacci numbers
$$
\Phi_{n,k}^{(s)} = 5^{-s}(F_k)^{-2s - 1}\sum_{j = 0}^{s}(-1)^{(k - 1)j}{n + 2s \choose j}{n + s - 1 - j \choose n - 1}F_{k(n + 2s - 2j)};
$$
\item via the recurrence relation
\begin{itemize}
\item define the linear form $(1 - L_kz + (-1)^kz^2)^{s + 1} = \Lambda(1,z,z^2,\ldots,z^{2s + 2})$;
\item then
$$
\Lambda(x_{n + 2s + 2}, x_{n + 2s + 1}, \ldots, x_n) = 0, \ x_j = \Phi_{j,k}^{(s)}, \ j = 1,\ldots, n + 2s + 1.
$$
\end{itemize}
\end{itemize}

\section{Further Identities}

Identities connecting Fibonacci number, Lucas nubmers, and the binomial coefficients have been known for a long time, and this community is replenished almost every day \cite{ref7,ref15,ref16,ref17,ref18,ref19}. In this section, consequences of the identities obtained below will be presented. Of course, most of these identities are known (or are easily derived from known ones), but, nevertheless, we hope that some identities are new.

From formulas \eqref{eq4} we obtain
$$
\sum_{j = 0}^{n - 1}F_{j + 1}F_{n - j} = \frac{1}{5}(nF_{n + 2} + (n + 2)F_n),
$$
$$
\sum_{j = 0}^{n - 1}F_{j + 1}\left((n - j)F_{n + 2 - j} + (n + 2 - j)F_{n - j}\right) = 
$$
$$
\frac{1}{5}\left(\frac{n(n + 1)}{2}F_{n + 4} + n(n + 4)F_{n + 2} + \frac{(n + 3)(n + 4)}{2}F_n\right),
$$
and then
$$
\sum_{j = 0}^{n - 1}(n - j)F_{j + 1}L_{n + 1 - j} = \frac{1}{10}\left(n(5n + 7)F_{n + 2} + 2(n + 2)F_n\right).
$$

Using the formula from example 18.7 \cite[p.~222]{ref7}, we obtain
$$
\sum_{j = 1}^{\left[n/2\right]}j{n - j \choose j} = \frac{1}{5}\left((n - 1)F_{n + 1} + (n + 1)F_{n - 1}\right) = \frac{1}{5}(nL_n - F_n).
$$

From \eqref{eq18}, taking into account that $F_1^{(s)} = 1$, $F_2^{(s)} = s + 1$, $F_3^{(s)} = \frac{1}{2}(s + 1)(s + 4)$, $F_4^{(s)} = \frac{1}{6}(s + 1)(s + 2)(s + 9)$, and $F_5^{(s)} = \frac{1}{24}(s + 1)(s + 2)(s + 4)(s + 15)$, it follows that
$$
\frac{1}{5^s}\sum_{j = 0}^s{1 + 2s \choose j}F_{1 + 2s - 2j} = 1,
$$
$$
\frac{1}{5^s}\sum_{j = 0}^s{2 + 2s \choose j}(1 + s - j)F_{2 + 2s - 2j} = s + 1,
$$
$$
\frac{1}{5^s}\sum_{j = 0}^s{3 + 2s \choose j}{2 + s - j \choose 2}F_{3 + 2s - 2j} = \frac{1}{2}(s + 1)(s + 4),
$$
$$
\frac{1}{5^s}\sum_{j = 0}^s{4 + 2s \choose j}{3 + s - j \choose 3}F_{4 + 2s - 2j} = \frac{1}{6}(s + 1)(s + 2)(s + 9),
$$
$$
\frac{1}{5^s}\sum_{j = 0}^s{5 + 2s \choose j}{4 + s - j \choose 4}F_{5 + 2s - 2j} = \frac{1}{24}(s + 1)(s + 2)(s + 4)(s + 15),
$$
etc. for all natural numbers $s$.

From the definition of numbers $\Phi_{n,k}^{(s)}$ (formula (5)) we obtain:
$$
\sum_{j = 0}^{n - 1}F_{(j + 1)k}F_{(n - j)k} = \frac{1}{5F_k}\left(nF_{(n + 2)k} + (-1)^{k - 1}(n + 2)F_{nk}\right),
$$
$$
\sum_{j = 0}^{n - 1}F_{(j + 1)k}\left((n - j)F_{(n - j + 2)k} + (-1)^{k - 1}(n - j + 2)F_{(n - j)k}\right) =
$$
$$
\frac{1}{5F_k}\left(\frac{n(n + 1)}{2}F_{(n + 4)k} + (-1)^{k - 1}n(n + 4)F_{(n + 2)k} + \frac{(n + 3)(n + 4)}{2}F_{nk}\right),
$$
and then
$$
\sum_{j = 0}^{n - 1}(n - j)F_{(j + 1)k}\left(F_{(n - j + 2)k} + (-1)^{k - 1}F_{(n - j)k}\right) =
$$
$$
\frac{1}{5F_k}\left(\frac{n(n + 1)}{2}F_{(n + 4)k} + (-1)^{k - 1}n(n + 2)F_{(n + 2)k} + \frac{1}{2}(n^2 + 3n + 4)F_{nk}\right).
$$

Let us now take into account that $\Phi_{1,k}^{(s)} = 1$, $\Phi_{2,k}^{(s)} = (s + 1)L_k$, and $\Phi_{3,k}^{(s)} = \frac{1}{2}(s + 1)(s + 2)L_{2k} + (s + 1)^2(-1)^k$. Then, $\Phi_{4,k}^{(s)} = \frac{1}{2}(s + 1)(s + 2)(s + 3)(L_k)^3 - (s + 1)(s + 2)(-1)^kL_k$ by virtue of equalities
$$
\frac{1}{2^ss!}U_{2 + s}^{(s)}(z) = 2(s + 1)(s + 2)z^2 - (s + 1),
$$
$$
\frac{1}{2^ss!}U_{3 + s}^{(s)}(z) = \frac{4}{3}(s + 1)(s + 2)(s + 3)z^3 - 2(s + 1)(s + 2)z.
$$

Then,
\begin{equation}\label{eq23}
5^{-s}(F_k)^{-2s - 1}\sum_{j = 0}^{s}(-1)^{(k - 1)j}{1 + 2s \choose j}F_{k(1 + 2s - 2j)} = 1, 
\end{equation}
or
$$
(F_k)^{2s + 1} = 5^{-s}\sum_{j = 0}^{s}(-1)^{(k - 1)j}{1 + 2s \choose j}F_{k(1 + 2s - 2j)},
$$
\begin{equation}\label{eq24}
5^{-s}(F_k)^{-2s - 1}\sum_{j = 0}^{s}(-1)^{(k - 1)j}{2 + 2s \choose j}(1 + s - j)F_{k(2 + 2s - 2j)} = (s + 1)L_k, 
\end{equation}
or
$$
(F_k)^{2s} = \frac{5^{-s}}{(s + 1)F_{2k}}\sum_{j = 0}^{s}(-1)^{(k - 1)j}{2 + 2s \choose j}(1 + s - j)F_{2k(1 + s - j)},
$$
\begin{align}\label{eq25}
5^{-s}(F_k)^{-2s - 1}\sum_{j = 0}^{s}(-1)^{(k - 1)j}{3 + 2s \choose j}{2 + s - j \choose 2}F_{k(3 + 2s - 2j)} = \nonumber \\
\frac{1}{2}(s + 1)(s + 2)L_{2k} + (s + 1)^2(-1)^k,
\end{align}
\begin{align}\label{eq26}
5^{-s}(F_k)^{-2s - 1}\sum_{j = 0}^{s}(-1)^{(k - 1)j}{4 + 2s \choose j}{3 + s - j \choose 3}F_{k(4 + 2s - 2j)} = \nonumber \\
\frac{1}{6}(s + 1)(s + 2)(s + 3)(L_k)^3 - (s + 1)(s + 2)(-1)^kL_k. 
\end{align}

For small values of $k$ or $s$ the formulas can be obtained, which, as a rule, have already been discovered earlier (most likely by other methods). For example from \eqref{eq23} and \eqref{eq24}, by taking into account that $L_k = F_{2k}/F_k$ (for $s = 1,2,3$), we get 
$$
F_{2k}(F_k)^2 = \frac{1}{5}\left(F_{4k} + 2(-1)^{k - 1}F_{2k}\right),
$$
$$
F_{2k}(F_k)^4 = \frac{1}{25}\left(F_{6k} + 4(-1)^{k - 1}F_{4k} + 5F_{2k}\right),
$$
$$
F_{2k}(F_k)^6 = \frac{1}{125}\left(F_{8k} + 6(-1)^{k - 1}F_{6k} + 14F_{4k} + 14(-1)^{k - 1}F_{2k}\right),
$$
$$
L_k = \frac{F_{4k} - 2(-1)^kF_{2k}}{F_{3k} - 3(-1)^kF_k} = \frac{F_{6k} - 4(-1)^kF_{4k} + 5F_{2k}}{F_{5k} - 5(-1)^kF_{3k} + 10F_k},
$$
etc.

From \eqref{eq25} for $s = 0, 1, 2, 3$:
$$
\frac{F_{3k}}{F_k} = L_{2k} + (-1)^k,
$$
$$
\frac{3F_{5k} + 5(-1)^{k - 1}F_{3k}}{5(F_k)^3} = 3L_{2k} + 4(-1)^k,
$$
$$
\frac{2F_{7k} + 7(-1)^{k - 1}F_{5k} + 7F_{3k}}{25(F_k)^5} = 2L_{2k} + 3(-1)^k,
$$
$$
\frac{5F_{9k} + 27(-1)^{k - 1}F_{7k} + 54F_{5k} + 42(-1)^{k - 1}F_{3k}}{125(F_k)^7} = 5L_{2k} + 8(-1)^k,
$$
etc.

From \eqref{eq26} for $s = 0,1,2,3$:
$$
\frac{F_{4k}}{F_k} = (L_k)^3 - 2(-1)^kL_k,
$$
$$
\frac{2F_{6k} + 3(-1)^{k - 1}F_{4k}}{5(F_k)^3} = 2(L_k)^3 - 3(-1)^kL_k,
$$
$$
\frac{5F_{8k} + 16(-1)^{k - 1}F_{6k} + 14F_{4k}}{25(F_k)^5} = 5(L_k)^3 - 6(-1)^kL_k
$$
$$
\frac{F_{10k} + 5(-1)^{k - 1}F_{8k} + 9F_{6k} + 6(-1)^{k - 1}F_{4k}}{125(F_k)^7} = (L_k)^3 - (-1)^kL_k,
$$
etc.

Using formula \eqref{eq22}, one can also obtain further identities for different values of $n$, $k$, and $s$.

\section{Conclusion}

The role of Chebyshev polynomials in Mathematics can not be overstated: they are a powerful tool for solving a variety of problems in interpolation theory, approximation theory, numerical analysis, dynamical systems theory, number theory, etc. However, derivatives of Chebyshev polynomials, especially higher orders, appear considerably less frequently in the literature, primarily in studies of the general properties of orthogonal polynomials. Some publications can be cited where they play an independent role, for example \cite{ref20,ref21,ref22}. In this article, such polynomials are used as the basis for solving certain problems in number theory.

\section{Acknowledgment}

A. Stokolos is partially supported by the AMS-Simons Research Enhancement Grant \#367012 for PUI faculty.

\bibliography{seqbib}
\bibliographystyle{amsplain}

\end{document}